\documentclass[10pt]{amsart}
\usepackage{amsthm, amsfonts, amssymb, amsmath, amscd}
\usepackage[all]{xy}
\usepackage{lmodern}
\usepackage[margin=3.5cm]{geometry}

\newcommand{\To}{\rightarrow}

\newtheorem*{thm*}{Theorem}

\theoremstyle{definition}

\begin{document}
\title{The Mazur - Ulam theorem}
\author{Bogdan Nica}
\address{Mathematisches Institut, Georg-August-Universit\"at G\"ottingen, Germany}
\email{bogdan.nica@gmail.com}
\date{\today}

\maketitle
In this note, we give a short proof of the following classical theorem:

\begin{thm*}[Mazur - Ulam 1932]
Every bijective isometry between real normed spaces is affine.
\end{thm*}

Recall that a map $\alpha:X\To Y$ between real normed spaces is said to be \emph{affine} if it satisfies $\alpha((1-t)x_1+tx_2)=(1-t)\alpha(x_1)+t\alpha(x_2)$ for all $x_1,x_2\in X$ and $t\in [0,1]$. This definition turns out to be equivalent to the requirement that $\alpha$ is linear up to a translation, i.e., $\alpha-\alpha(0)$ is a linear map. In applications, the second characterization is often the useful one. For instance, in geometric group theory one is interested in isometric actions of groups on Banach spaces. The Mazur - Ulam theorem implies that every isometric group action on a real Banach space is given by a linear isometric action and a group cocycle with values in the given Banach space.

Our proof of the Mazur - Ulam theorem is inspired by V\"{a}is\"{a}l\"{a}'s proof \cite{Va03}, which in turn is based on ideas of Vogt \cite{Vogt}. See for instance Lax \cite[pp.49-51]{Lax} for an exposition of the original proof by Mazur and Ulam \cite{MU}. A common point to all these approaches is the use of reflections. Compared to the proof of Mazur and Ulam, or to that of V\"{a}is\"{a}l\"{a}, the proof we give below needs no preliminary setup. The reflection trick is all the more striking.

\begin{proof} Let $X$ be a real normed space, and fix two points $x_1,x_2\in X$. We check that $\alpha(\frac{x_1+x_2}{2})=\frac{\alpha(x_1)+\alpha(x_2)}{2}$ for every bijective isometry $\alpha$ defined on $X$. For any such $\alpha$, we let
\begin{eqnarray*}
\mathrm{def}(\alpha)=\Big\|\alpha\Big(\frac{x_1+x_2}{2}\Big)-\frac{\alpha(x_1)+\alpha(x_2)}{2}\Big\|
\end{eqnarray*} 
denote the possible ``affine defect''. Observe that we have a uniform bound on the defect:
\begin{displaymath}
\mathrm{def}(\alpha)\leq\frac{1}{2}\Big\|\alpha\Big(\frac{x_1+x_2}{2}\Big)-\alpha(x_1)\Big\| +\frac{1}{2}\Big\|\alpha\Big(\frac{x_1+x_2}{2}\Big)-\alpha(x_2)\Big\|=\frac{\|x_1-x_2\|}{2}
\end{displaymath}
Now here comes the magic: for every bijective isometry $\alpha$ defined on $X$ we define another bijective isometry $\alpha'$ on $X$ whose defect is twice the defect of $\alpha$. Indeed, the bijectivity of $\alpha$ allows us to consider $\alpha'=\alpha^{-1}\rho \alpha$ where $\rho$ is the reflection in $\frac{\alpha(x_1)+\alpha(x_2)}{2}$ in the target space of $\alpha$. Thus $\rho$ is given by $z\mapsto \alpha(x_1)+\alpha(x_2)-z$, and $\alpha'(x_1)=x_2$, $\alpha'(x_2)=x_1$. As $\alpha^{-1}$ is an isometry, we have:
\begin{eqnarray*}
\mathrm{def}(\alpha')&=&\bigg\|\alpha^{-1}\bigg(\alpha(x_1)+\alpha(x_2)-\alpha\Big(\frac{x_1+x_2}{2}\Big)\bigg)-\frac{x_1+x_2}{2}\bigg\|\\
&=&\Big\|\alpha(x_1)+\alpha(x_2)-\alpha\Big(\frac{x_1+x_2}{2}\Big)-\alpha\Big(\frac{x_1+x_2}{2}\Big)\Big\|=2\;\mathrm{def} (\alpha)
\end{eqnarray*}
If we had a bijective isometry defined on $X$ with positive affine defect, then by iteration we would obtain bijective isometries defined on $X$ with arbitrarily large affine defect, and that would contradict the uniform bound we established in the first place. We conclude that $\mathrm{def} (\alpha)=0$ for every bijective isometry $\alpha$ defined on $X$. 

Once we know that a bijective isometry $\alpha$ defined on $X$ satisfies $\alpha(\frac{x_1+x_2}{2})=\frac{\alpha(x_1)+\alpha(x_2)}{2}$ for all $x_1,x_2\in X$, a standard continuity argument shows that $\alpha$ is affine. 
\end{proof}

\noindent\textbf{Acknowledgment.} I thank Jussi V\"{a}is\"{a}l\"{a} for his interest in this note.


\begin{thebibliography}{00}
\bibitem{Lax} P. Lax: \emph{Functional analysis}, Wiley-Interscience Series in Pure and Applied Mathematics 2002
\bibitem{MU} S. Mazur, S. Ulam: \emph{Sur les transformations isom\'{e}triques d'espaces vectoriels norm\'{e}s}, C. R. Acad. Sci. Paris 194 (1932), 946--948
\bibitem{Va03} J. V\"{a}is\"{a}l\"{a}: \emph{A proof of the Mazur-Ulam theorem}, Amer. Math. Monthly 110 (2003), no. 7, 633--635
\bibitem{Vogt} A. Vogt: \emph{Maps which preserve equality of distance}, Studia Math. 45 (1973) 43--48.
\end{thebibliography}
\end{document}